\documentclass[12pt]{article}
\usepackage{amsfonts,amssymb,amsmath}
\usepackage[english]{babel}
\usepackage{a4wide}
\usepackage{amssymb}
\usepackage{amsthm}
\usepackage{pdfsync}

\usepackage{numbysec}

-\marginparwidth \oddsidemargin -4mm \evensidemargin \oddsidemargin

\newcommand{\commentout}[1]{}

\newcommand{\bbE}{\mathbb{E}}
\newcommand{\bbP}{\mathbb P}

\newcommand{\bbX}{X}

\newcommand{\bbV}{V}

\newcommand{\bbx}{x}

\newcommand{\la}{\lambda}

\newtheorem{thm}{Theorem}[section]

\newtheorem{proposition}[theorem]{Proposition}

\newcommand{\bes}{\begin{displaymath}}
\newcommand{\ees}{\end{displaymath}}
\newcommand{\be}{\begin{equation}}
\newcommand{\ee}{\end{equation}}
\newcommand{\ba}{\begin{eqnarray}}
\newcommand{\ea}{\end{eqnarray}}
\newcommand{\bas}{\begin{eqnarray*}}
\newcommand{\eas}{\end{eqnarray*}}

\newcommand{\B}{{\@Bbb B}}
\newcommand{\C}{{\@Bbb C}}

\newcommand{\F}{{\@Bbb F}}
\renewcommand{\P}{{\@Bbb P}}

\newcommand{\Q}{{\@Bbb Q}}
\newcommand{\bQ}{{\@Bbb Q}}
\newcommand{\N}{{\@Bbb N}}
\newcommand{\R}{{\@Bbb R}}

\newcommand{\W}{{\@Bbb W}}

\newcommand{\al}{\alpha}

\newcommand{\vv}{v}

\newcommand{\bze}{0}

\newcommand{\cA}{\@s A}
\newcommand{\cB}{\@s B}
\newcommand{\cC}{\@s C}
\newcommand{\cD}{\@s D}
\newcommand{\cE}{\@s E}
\newcommand{\cF}{\@s F}
\newcommand{\cG}{\@s G}
\newcommand{\cH}{\@s H}
\newcommand{\cI}{\@s I}
\newcommand{\cJ}{\@s J}
\newcommand{\cK}{\@s K}
\newcommand{\cL}{\@s L}
\newcommand{\cN}{\@s N}
\newcommand{\cM}{\@s M}
\newcommand{\cO}{\@s O}
\newcommand{\cP}{\@s P}
\newcommand{\cR}{\@s R}
\newcommand{\cS}{\@s S}
\newcommand{\cT}{\@s T}
\newcommand{\cV}{\@s V}
\newcommand{\cW}{\@s W}
\newcommand{\cX}{\@s X}
\newcommand{\cY}{\@s Y}
\newcommand{\cZ}{\@s Z}

\newcommand{\bma}{\@bm a}
\newcommand{\bmb}{\@bm b}
\newcommand{\bmc}{\@bm c}
\newcommand{\bmd}{\@bm d}

\newcommand{\bme}{\@bm e}
\newcommand{\bmf}{\@bm f}
\newcommand{\bmg}{\@bm g}
\newcommand{\bmh}{\@bm h}
\newcommand{\bmi}{\@bm i}
\newcommand{\bmj}{\@bm j}
\newcommand{\bmk}{\@bm k}
\newcommand{\bml}{\@bm l}
\newcommand{\bmm}{\@bm m}
\newcommand{\bmn}{\@bm n}
\newcommand{\bmo}{\@bm o}
\newcommand{\bmp}{\@bm p}
\newcommand{\bmq}{\@bm q}
\newcommand{\bmr}{\@bm r}
\newcommand{\bms}{\@bm s}
\newcommand{\bmt}{\@bm t}
\newcommand{\bmu}{\@bm u}
\newcommand{\bmw}{\@bm w}
\newcommand{\bmv}{\@bm v}
\newcommand{\bmx}{\@bm x}
\newcommand{\bx}{\@bm x}
\newcommand{\bmy}{\@bm y}
\newcommand{\bmz}{\@bm z}

\newcommand{\by}{\@bm y}
\newcommand{\bmzero}{\@bm 0}

\newcommand{\gA}{\@g A}
\newcommand{\gD}{\@g D}
\newcommand{\gJ}{\@g J}
\newcommand{\gF}{\@g F}
\newcommand{\gM}{\@g M}
\newcommand{\gR}{\@g R}

\newcommand{\bR}{R}
\newcommand{\mbR}{\mathbb{R}}
\newcommand{\bbC}{\mathbb{C}}

\title{Passive tracer in non-Markovian, Gaussian velocity field}

\begin{document}
\numberbysection

\author{Tymoteusz Chojecki\thanks{Institute of Mathematics,  UMCS,
pl. Marii Curie-Sk\l odowskiej 1, 20-031, Lublin}}

\maketitle
\begin{abstract}
We consider the trajectory of a tracer  that is  the solution of an
ordinary differential equation $\dot\bbX(t)=\bbV(t, \bbX(t)),\
X(0)=0$, with the right hand side, that is a stationary, zero-mean,
Gaussian vector field with incompressible realizations. It is known,
see \cite{komorowski-fannjiang,carmona-xu,Ksiazka}, that
$\bbX(t)/\sqrt{t}$ converges in law, as $t\to+\infty$, to a normal, zero mean
vector, provided that the field $V(t,x)$ is Markovian and has the
 spectral gap property. We wish to extend this result to the case
 when the field is not Markovian and its covariance matrix is given by a completely monotone Bernstein function.
\end{abstract}

\section{Introduction and some assumptions}\label{wstep}
In this paper we would like to show the central limit theorem for a
passive tracer model, when the velocity field is non-Markovian but
Gaussian and exponentially mixing in time.

Passive tracer model is given by the following equation,
\begin{equation}\label{lila4}
\left\{\begin{array}{ll}
&\dfrac{dX(t)}{dt}=V\left(t,X(t)\right),\quad t>0,\\
&\\
 &X(0)=0,\end{array} \right.
\end{equation}
where $\bbV:\mbR^{1+d}\times\Omega\to\mbR^d$ is a real,
$d-$dimensional, incompressible i.e.
$\sum_{p=1}^d\partial_{x_p}V_p(t,x)\equiv0,$ zero mean, Gaussian
random vector field over a probability space
$(\Omega,\mathcal{F},\bbP)$.

Some basic problems concerning the asymptotic behavior of the tracer
are: the  law of large numbers (LLN) i.e. whether $X(t)/t$ converges
to a constant vector $v_*$ (called the  {\em Stokes drift}), as
$t\to+\infty$ and the central limit theorem (CLT), i.e. whether $
(X(t)-\vv_* t)/\sqrt{t}$ is convergent in law to a normal vector
$N(0,\kappa)$. The covariance matrix
$\kappa=[\kappa_{ij}]_{i,j=1,\ldots d}$ is called {\em turbulent
diffusivity} of the tracer.

It is expected that both the LLN, with $v_*=0$, and the CLT for the tracer trajectory hold when the velocity field is zero mean, Gaussian, incompressible  and its
 covariance matrix $R(t,x)=[R_{pq}(t,x)]_{p,q=1,\dots,d}$, given by
 $$
R_{pq}(t,x):=\bbE[V_p(t,x)V_q(0,0)],\quad p,q=1,\dots,d,\
(t,x)\in\mbR^{1+d},
$$
   exponentially decays in
time, i.e. there exists $C>0$ such, that
\begin{equation}\label{chciec}
\sum_{p,q=1}^d|R_{pq}(t,\bbx)|\leq Ce^{-|t|/C},\ \textrm{for all }
(t,x)\in\mbR^{1+d}.
\end{equation}

The CLT has been established in \cite{Prof}, in the case of
$T-$dependent fields, i.e. those for which exists $T>0$ such that
$\bR(t,x)=0,\ |t|>T,\ x\in\mbR^d$.

In case when the vector field $ \bbV(t,x)$ is Markovian (not
necessarily Gaussian) and satisfies the spectral gap condition, the
CLT has been established in  \cite{komorowski-fannjiang}, Theorem A,
see also \cite{Ksiazka,carmona-xu,koralov}.
 In the Gaussian case when the covariance matrix is of the form
 \begin{equation}\label{KowariancjaMarkowa}
R_{pq}(t,\bbx)=\int_{\mbR^d}e^{i\bbx\cdot\xi-\gamma(\xi)|t|}\hat
R_{pq}(d\xi),\quad p,q=1,\dots,d,\ (t,x)\in\mbR^{1+d},
\end{equation}
 where both $\gamma(\cdot)$ and non-negative Hermitian matrix valued measure $\hat R(\cdot)=[\hat R_{pq}(\cdot)]$ are even (because the field is real), then the field is Markovian.
  It can be shown, see Chapter 12 of \cite{Ksiazka}, that the spectral gap condition holds, provided  there exists $\gamma_0>0$ such that
\begin{equation}\label{gammaWar}
\gamma(\xi)\ge\gamma_0,\quad\xi\in\mbR^d.
\end{equation}

We will consider fields for which the exponential factor is replaced
by a function $h:[0,+\infty)\times\mbR^d\to\mbR$ i.e. fields with the
covariance defined as follows
\begin{equation} \label{Kowariancja}
 R_{pq}(t,\bbx)=\int_{\mbR^d}e^{i\bbx\cdot\xi}h(|t|,\xi)\hat R_{pq}(\xi)d\xi,\qquad p,q=1,\dots,d,\ (t,x)\in\mbR^{1+d},
\end{equation}
where $\hat R_{pq}(\xi)$ is a density of $\hat R_{pq}(\cdot)$.

 We show in Proposition \ref{dodatnioh} that the function $h$ is
non-negative definite iff $R(t,x)=[R_{pq}(t,x)]$ is non-negative
definite. Therefore, the largest (in the sense of inclusion) set of
functions $h$ in \eqref{Kowariancja} which can be examined is the
class of non-negative definite functions. We study a smaller family of
functions, namely we assume that $h$ is \emph{completely monotone in
the sense of Bernstein}, see \eqref{hfunkcja1}.

Let us denote by $\hat r$ the \emph{power-energy spectrum}. It is a
scalar non-negative, integrable function given by formula
\begin{equation}\label{power-energy}
\hat r(\xi):=\textrm{tr}\hat R(\xi),\quad \xi\in\mbR^d,
\end{equation}
where tr is the trace. Let $
\frak{r}(d\xi):=\hat r(\xi)d\xi.
$

The main result of the paper, see Theorem \ref{Glowne}, is the CLT
for the trajectory of a tracer moving in a field whose covariance
matrices are given by \eqref{Kowariancja}, where the function $h$ is
\emph{completely monotone} i.e. $h\in
C^{\infty}(0,+\infty)\cap C[0,+\infty),$ $(-1)^nh^{(n)}(t,\xi)\ge
0,\ t>0,\ {\frak r}\ \textrm{a.e.}\ \xi,\ n=0,1,\dots$, and satisfies
\eqref{miaranos}. From the Bernstein Theorem (\cite{Lax} Theorem 3., p.
138) we know that
\begin{equation}\label{hfunkcja1}
h(t,\xi)= \int_{0}^{+\infty}e^{-\lambda t}\mu(\xi,d\lambda),\qquad
t\in[0,+\infty),\ {\frak r}\ \textrm{a.e.}\ \xi,
\end{equation}
where $\mu(\xi,\cdot)$ is a non-negative, finite Borel measure on
$[0,+\infty)$ for ${\frak r}\ \textrm{a.e.}\ \xi$. For example from
\cite{Vondracek}, Lemma 4.5, we know, that all completely monotone
functions are non-negative definite. We assume that there exists
$\lambda_0>0$ such that
\begin{equation}\label{miaranos}
\textrm{supp } \mu(\xi,\cdot)\subset[\lambda_0,+\infty),\quad {\frak r}\ \textrm{a.e.}\ \xi.
\end{equation}
 Observe that this assumption implies
\eqref{chciec}.

In Section \ref{exy} we show (see \eqref{przyklad}) an example of a
covariance matrix which is of the form \eqref{Kowariancja} but
not of the form \eqref{KowariancjaMarkowa}.

We prove Theorem \ref{Glowne} in Section \ref{dow} by embedding the
field $V(t,x)$ into a larger space where we add one dimension and
one argument i.e a space of $d+1$ dimensional  fields $\tilde
V(t,x,y)$. We define a field $\tilde V(t,x,y)$ in such a way that
the field $(\tilde V(t,x,0))$ has the same distribution as the field
$(V(t,x),0)$. The process $\tilde V(t,\cdot,\cdot)$ has the Markov
property in appropriate function space. This process also has the
spectral gap property. At the end we use Theorem 12.13 from
\cite{Ksiazka}.

In Section \ref{pomocniczedowody} we show the proof of Proposition
\ref{dodatnioh}.

\section{Preliminaries and the statement of the main result}
First, we present some assumptions on matrix valued function
$
\hat R(\xi)=[\hat R_{pq}(\xi)]_{p,q=1,\dots,d},
$
 which guarantee
that  function $R(t,x)$ (defined in \eqref{Kowariancja}) is non-negative definite:%
%
\begin{equation}\label{zalb}
\hat R_{pq}(\xi)=\hat R_{qp}^*(\xi),\quad p,q=1,\dots,d,\ {\frak r}\
\textrm{a.e.}\ \xi,
\end{equation}
\begin{equation}\label{niuejemnahat}
\hat R(\xi)\eta\cdot\eta\geq0,\quad \eta\in\bbC^d,\ {\frak r}\
\textrm{a.e.}\ \xi.
\end{equation}


Now we present the result which explains why we need to deal with
functions $h$, which are non-negative definite in $t$ for ${\frak
r}$ $a.e$ $\xi$.
\begin{proposition}\label{dodatnioh}
For any $N\ge1$, $\al_1,\ldots,\al_N\in\mathbb C$ and
$t_1,\ldots,t_N\in\mbR$ we have
\begin{equation}
\label{pos-def} \sum_{p,q=1}^Nh(|t_p-t_q|,\xi)\al_p\bar
\al_q\ge0,\quad \textrm{for}\ {\frak r}\ \textrm{a.e.}\ \xi,
\end{equation}
iff the matrix valued function $R$ given by \eqref{Kowariancja} is
non-negative definite.
\end{proposition}{\hspace*{-2em}}
The proof of this Proposition is presented in Section \ref{pomocniczedowody}.

 Since $\bbE|V(0,0)|^2<+\infty,$ we have
\begin{equation}\label{zalc}
\int_{\mbR^d}h(0,\xi)\hat r(\xi)d\xi<+\infty.
\end{equation}
Recall that we have assumed that $h$ is of the form \eqref{hfunkcja1}. To
fulfill the assumption \eqref{zalc} we require that
\begin{equation}\label{esssup}
\mathop{\textrm{esssup }}_{\xi}\mu(\xi,[0,+\infty))<+\infty.
\end{equation}
It implies that $|h|$ is bounded. Above esssup is the essential
supremum with respect to ${\frak r}$. Then \eqref{zalc} is implied by
\begin{equation}\label{zalcc}
\int_{\mbR^d}\hat r(\xi)d\xi<+\infty.
\end{equation}

Assumptions \eqref{zalb}--\eqref{zalcc} imply that the matrix valued function
$R(\cdot,\cdot)$ is non-negative definite. From \cite{rozanov}, Section I.3,
we know that there exits a unique (in the sense of law)
stationary Gaussian random vector field $V(t,x)$ such that $R(t,x)$
is its covariance matrix.

We want to deal with a real field, so we need to assume, see
\cite{Rozanow2} Theorem 4.2., p. 18, that
\begin{itemize}
\item For ${\frak r}$ a.e. $\xi$ we have
\begin{equation}\label{zald}
\hat R_{pq}(\xi)=\hat R_{qp}(-\xi),\quad p,q=1,\dots,d.
\end{equation}
\item For any $t\in\mbR$
\begin{equation}\label{parzxi}
h(|t|,\xi)=h(|t|,-\xi),\quad {\frak r}\ \textrm{a.e.}\ \xi.
\end{equation}
\end{itemize}

 Assumption \eqref{parzxi} implies (see \eqref{hfunkcja1})
\begin{equation}\label{parzximu}
\mu(\xi,\cdot)=\mu(-\xi,\cdot),\quad{\frak r}\ \textrm{a.e.}\ \xi.
\end{equation}

To sum up, in this paper we consider the model in which a
$d$-dimensional random vector field  $\bbV(t,x)$ has the covariance
matrix given by \eqref{Kowariancja}, where $\hat R(\xi)$ and $\
h(|t|,\xi)$ satisfy assumptions \eqref{zalb}--\eqref{parzximu}.

To be able to solve the equation \eqref{lila4} we need to assume an
appropriate regularity of the field. Namely, we assume that there
exists the second derivative in $x$ of the field $\bbV(t,x)$ and it
is continuous. This implies (see \eqref{Kowariancja})
\begin{equation}
\int_{\mbR^d}(1+|\xi|^4)\hat r(\xi)d\xi<+\infty.\label{IIPoch}
\end{equation}

Moreover, we assume that the field $V(t,x)$ is incompressible, i.e.
$$
\nabla_x\cdot V(t,x)=\sum_{p=1}^d\partial_{x_p}V_p(t,x)\equiv0,\quad
(t,x)\in\mbR^{1+d}.
$$
 We
can rewrite this using \eqref{Kowariancja} as follows
\begin{equation}\label{A2}
\sum_{p=1}^d\xi_p\hat R_{pq}(\xi)=0,\quad q=1,\dots,d,\ {\frak r}\
\textrm{a.e.}\ \xi.
\end{equation}
From the above equality we obtain
\begin{equation}\label{A3}
\sum_{p=1}^d\partial_{x_p}R_{pq}(t,x)=0,\quad (t,x)\in\mbR^{1+d},\
q=1,\dots,d.
\end{equation}


Formula \eqref{Kowariancja} and condition \eqref{IIPoch} imply in
particular (see Theorem 3.1 from \cite{watanabe}) that there exists
a version of $V(t,x)$ with the realizations which grows slower than
linearly in $(t,x)$ a.s., therefore the equation \eqref{lila4} can
be solved globally in $t$, and the process $X(t)$ is defined for all
$t\in[0,+\infty)$.

Now we can present our main result of this paper.

\begin{thm}\label{Glowne}
Assume that $V(t,x)$ is a zero-mean Gaussian random vector field
whose covariance matrix $R(t,x)=[R_{pq}(t,x)]_{p,q=1,\dots,d}$ is
given by \eqref{Kowariancja}. In addition suppose that \eqref{zalb}--\eqref{A2} hold. Then
the random variables $X(t)/\sqrt{t}$, where $X(t)$ is defined
  in $\eqref{lila4}$, converge weakly to the normal vector $N(\bze,\kappa)$, where
   the limit covariance matrix $\kappa=[\kappa_{pq}]$ satisfies:
$$
\kappa_{pq}=\lim_{t\to+\infty}\frac{1}{t}\bbE[X_p(t)X_q(t)]\quad
p,q=1,\dots,d.
$$
\end{thm}
\section{Proof of Theorem \ref{Glowne}} \label{dow}
%


Consider a $d+1-$dimensional field  $\widetilde
\bbV:\mbR^{2+d}\times\Omega\to\mbR^{1+d}$ with the covariance matrix
\begin{equation}
\label{SakoTRG} \widetilde R_{pq}(t,x,y)=\int_{\mbR}\int_{\mbR^d}
e^{-|\lambda t|}e^{i\xi\cdot x+i\lambda y}\hat
R_{pq}(\xi)m(d\xi,d\lambda),\quad t,y\in\mbR,x\in\mbR^d,
\end{equation}
 for $p,q=1,\ldots,d$, $\widetilde R_{pq}=0,$ if $p$ or $q=d+1$.
Above $m(d\xi,d\lambda):=\tilde \mu(\xi,d\lambda)d\xi$ would be the
measure given by
\begin{equation}\label{R11L}
 m(A\times
B)=\int_Ad\xi\left\{\int_B\tilde\mu(\xi,d\lambda)\right\},\quad A\in
{\cal B}(\mbR^d),B\in{\cal B}(\mbR),
\end{equation}
where $\tilde\mu$ is defined as follows
$$
\tilde\mu(\xi,d\lambda):=\frac{1}{2}\textbf{1}_{(0,+\infty)}(\lambda)\mu(\xi,d\lambda)+
\frac{1}{2}\textbf{1}_{(-\infty,0]}(\lambda)\mu(\xi,-d\lambda),
$$
recall that $\mu$ is given by \eqref{hfunkcja1}. The covariance matrix
in full Fourier transform form is given by
$$
\widetilde
R_{pq}(t,x,y)=\frac{1}{\pi}\int_{\mbR}\int_{\mbR}\int_{\mbR^d}\frac{\lambda}{\lambda^2+\tau^2}e^{i\tau
t+ix\cdot\xi+i\lambda y}\hat R_{pq}(\xi)\tilde\mu(\xi,d\lambda)d\tau
d\xi.
$$
 This equality together
with the assumptions \eqref{zalb}-\eqref{A2} imply that $\widetilde
R_{pq}(t,x,y)$ is the covariance matrix of some zero mean,
Gaussian, stationary in time and space, random, real valued vector field
$\widetilde V(t,x,y)$ (\cite{rozanov}, Section I.3). Observe that
for $y=0$  we have
\begin{align}\label{rowkow}\begin{aligned}
&\widetilde R_{pq}(t,x,0)=\int_{\mbR}\int_{\mbR^d} e^{-|\lambda
t|}e^{i\xi\cdot x}\tilde\mu(\xi,d\lambda)\hat
R_{pq}(\xi)d\xi=\frac{1}{2}\int_{\mbR^d}\int_0^{+\infty}e^{-\lambda|t|}e^{i\xi\cdot
x}\mu(\xi,d\lambda)\hat R_{pq}(\xi)d\xi\\
&+\frac{1}{2}\int_{\mbR^d}\int_{-\infty}^{0}e^{\lambda|t|}e^{i\xi\cdot
x}\mu(\xi,-d\lambda)\hat R_{pq}(\xi)d\xi=R_{pq}(t,x), \quad
p,q=1,\dots,d,\ (t,x)\in\mbR^{1+d},
\end{aligned}\end{align}
where matrix $R_{pq}(t,x)$ is given by \eqref{Kowariancja}. Observe
that $\widetilde V_{d+1}(t,x,y)\equiv0$, because $\widetilde
R_{d+1,d+1}\equiv0$. So the field satisfies
\begin{equation}\label{pola}
 \left\{\widetilde\bbV(t,x,0),\ (t,x)\in\mbR^{1+d}\right\}\stackrel{d}{=}\left\{\left(\bbV(t,x),0\right),\ (t,x)\in\mbR^{1+d}\right\},
\end{equation}
where $\stackrel{d}{=}$ denotes equality in the law.
 Let $\tilde r$ be a measure
$
\tilde r(d\xi,d\lambda):=\tilde\mu(\xi,d\lambda)\hat r(\xi)d\xi.
$
From \eqref{miaranos} we know that $|\lambda|\geq\lambda_0,\
\tilde r$ a.e $\xi$.

Let $\rho>(d+1)/2$.
 Denote by $\mathcal{E}$ the Hilbert space
that is the completion of the space of functions
$v=(v_1,\dots,v_{d+1}):\mbR^{d+1}\to\mbR^{d+1}$ with components in
$C_c^{\infty}(\mbR^{d+1})$ (the space of infinitely differentiable
compactly supported functions), under the norm
$$
\|v\|^2_{\mathcal{E}}:=\int_{\mbR^{d+1}}\left(|v(z)|^2+|\nabla
v(z)|^2\right)(1+|z|^2)^{-\rho}dz.
$$
Here
$$
|v(z)|^2:=\sum_{i=1}^{d+1}v_i^2(z)\qquad\textrm{ and }\qquad|\nabla
v(z)|^2:=\sum_{i,j=1}^{d+1}\left(\partial_{z_j}v_i\right)^2(z),\quad
z\in\mbR^{d+1}.
$$
 Consider the
process
$$
\tilde V_t:=\tilde V(t,\cdot,\cdot),\ t\geq0.
$$
 From \cite{Ksiazka} Chapter 12, we know that the process takes
values in space $\mathcal{E}$. This space contains the realization
of the field $\tilde V(0,\cdot,\cdot)$. This is a separable Hilbert space
which turns out to be a subspace of the space of $C^1$ maps
$v:\mbR^{1+d}\to\mbR^{1+d}$. It can be shown that process
$\{\tilde V_t, t\ge0\}$ is stationary. On the space $\mathcal{E}$ we
introduce a measure $\pi$ as a distribution of $\tilde
V(0,\cdot,\cdot)$. Process $\{\tilde
V_t, t\ge0\}$ has the Markov property with the corresponding transition semigroup
$(P_t)_{t\ge0}$ i.e.
\begin{equation}\label{Markow2}
\bbE\left[F(\tilde V_{t+h})|\mathcal{V}_t\right]=P_hF(\tilde
V_t),\quad h\geq0,t\geq0,\ F\in\mathcal{B}_b(\mathcal{E}),
\end{equation}
where $\mathcal{V}_t:=\sigma(\tilde V_s,\ s\le t)$ is a natural
filtration of the process $\{\tilde V_t, t\ge0\}$. $(P_t)_{t\ge0}$ is a semigroup of Markov contractions on $L^2(\pi)$.
 Its generator $L:D(L)\to L^2(\pi)$ has the spectral gap
 property (see \cite{Ksiazka}, Corollary 12.15) i.e.
$$
-\langle LF,F\rangle_{L^2(\pi)}\ge \gamma_0\|F\|_{L^2(\pi)}^2,\quad
\textrm{for } F\in D(L),\textrm{ such that }\int_{{\cal E}}Fd\pi=0,
$$
where $\gamma_0$ was introduced in \eqref{gammaWar}.

We consider the following equation
\begin{equation}
 \label{lila3} \left\{\begin{array}{lllll}
\dfrac{d \tilde X_p(t)}{dt}=\widetilde \bbV_p(t,\tilde X(t),Y(t)),\quad p=1,\dots,d,\ t\ge0,&\\
&\\
\dfrac{dY(t)}{dt}=\widetilde\bbV_{d+1}(t,\tilde X(t),Y(t)),&\\ \tilde X(0)=0,&\\
Y(0)=0.
\end{array} \right.
\end{equation}
The existence and uniqueness of the solution \eqref{lila3} comes
from the form of the covariance of the field $\tilde V(\cdot)$.
Because $\widetilde V_{d+1}(t,x,y)\equiv 0$ we have $Y(t)=0$ for
$t\geq0$ . This and \eqref{pola} imply  $\{\tilde
X(t),t\ge0\}\stackrel{d}{=}\{X(t),t\ge0\}$, where $X(t)$ was defined
in \eqref{lila4}. So we will omit writing
 the symbol tilde over $X(t)$.


Now we introduce the so called environment process
$$
U_t:=\widetilde V(t,X(t)+\cdot,Y(t)+\cdot)\in\mathcal{E},\quad
 \bbP\ \textrm{a.e.}\ t\ge0.
$$
Observe that
$$
X(t)=\int_0^tF(U_s)ds,\quad t\ge0,
$$
where $F:{\cal E}\to\mbR^d$ is given by
$F(\omega)=(\omega_1(0,0),\ldots,\omega_d(0,0)),\ \omega\in{\cal
E}$. Process $\{U_t,t\ge0\}$ is also Markovian, see \cite{Ksiazka}
Proposition 12.19 (i). Observe that
$$
\nabla\cdot \widetilde V(t,x,y)=\sum_{i=1}^d\partial_{x_i}\widetilde
V_i(t,x,y)+\partial_y\widetilde
V_{d+1}(t,x,y)=\sum_{i=1}^d\partial_{x_i}\widetilde
V_i(t,x,y),\quad t,y\in\mbR,x\in\mbR^d.
$$
Let $Z(t,x,y):=\sum_{i=1}^d\partial_{x_i}\widetilde
V_i(t,x,y)$. Observe that
$$
\bbE Z^2(t,x,y)=\bbE
Z^2(0,0,0)=\sum_{i,j=1}^d\partial^2_{x_ix_j}\widetilde
R_{ij}(0,0,0)=\sum_{j=1}^d\partial_{x_j}\left\{\sum_{i=1}^d\partial_{x_i}R_{ij}(0,0)\right\}=0,
$$
where the last two equalities above come respectively from
\eqref{rowkow} and \eqref{A3}. Since $\nabla\cdot\widetilde V=0$
we know that $\pi$ is an invariant measure for the process
$\{U_t,t\geq0\}$, see \cite{Ksiazka} Proposition 12.19 (ii).
Moreover, the transition semigroup $(R_t)_{t\ge0}$ for
$\{U_t,t\ge0\}$, on $L^2(\pi)$,
with generator ${\cal L}:D({\cal L})\to L^2(\pi)$, has the spectral gap property
(\cite{Ksiazka} Corollary 12.22) i.e.
$$
-\langle {\cal L}F,F\rangle_{L^2(\pi)}\ge
\gamma_0\|F\|_{L^2(\pi)}^2,\quad \textrm{for } F\in D({\cal
L}),\textrm{ such that }\int_{{\cal E}}Fd\pi=0.
$$

The CLT for the process $\{X(t),t\ge0\}$ (which is
defined in \eqref{lila4}) comes directly from Theorem 12.13 of
\cite{Ksiazka}. \\
$\Box$

\subsection{Example of non-Markovian field which satisfies the assumption of Theorem \ref{Glowne}}\label{exy}

Observe that Theorem \ref{Glowne} can be applied to a field whose
covariance matrix is given by
\begin{equation}\label{przyklad}
R_{pq}(t,x)=\int_{\mbR^d}e^{ix\cdot\xi-\gamma(\xi)|t|-\tilde\gamma(\xi)|t|^\alpha}\hat
R_{pq}(d\xi),\ p,q=1,\dots,d,\ 0<\alpha\leq1,\ (t,x)\in\mbR^{1+d},
\end{equation}
where functions $\tilde\gamma,\gamma:\mbR^d\to\mbR$ satisfy
\eqref{gammaWar}.
 Observe that the function
\begin{equation}\label{talfa}
h(t,\xi)=\int_0^{+\infty}e^{-t^{\al}\la}\mu(\xi,d\la),\quad t\ge0,\ {\frak r}\
\textrm{a.e.} \xi,
\end{equation}
where $\al\in(0,1]$ is completely monotone in $t$. Indeed, we know
that
$$
H(t,\xi)=\int_0^{+\infty}e^{-t\la}\mu(\xi,d\la),\quad
t\ge0,\ {\frak r}\ \textrm{a.e.}\ \xi
$$
is completely monotone. However, $\phi(t):=t^{\al}$ can be
represented by an integral from $\psi(t):=\al t^{\al-1}$, which is
completely monotone on $(0,+\infty)$, if $\al\in(0,1]$. Therefore,
from Theorem 8 of \cite{schoenberg} we have that
$h(t,\xi)=H(\phi(t),\xi) $ is completely monotone on $[0,+\infty)$.

Now we will check condition \eqref{miaranos}. Observe that we can
write for $\mu(\xi,d\lambda):=\delta_{\tilde\gamma(\xi)}(d\lambda),\
t\geq0$,
$$
\int_0^{+\infty}e^{-t^\alpha\lambda}\mu(\xi,d\lambda)=e^{-\tilde\gamma(\xi)t^\alpha}=
\int_0^{+\infty}e^{-\tilde\gamma(\xi)\lambda
t}g(\lambda)d\lambda,\quad t\ge0,\  {\frak r}\ \textrm{a.e.}\ \xi,
$$
where  $g$ is a smooth density of $T(1)$, where $\{T(t),t \geq0\}$
is a $\alpha$-stable subordinator (\cite{Zolotarjev}, p.201, formula
(2.10.7)). Observe that for $t\geq0$
\begin{align}
&e^{-\gamma(\xi)t-\tilde\gamma(\xi)t^\alpha}=\int_0^{+\infty}e^{-(\tilde\gamma(\xi)\lambda+\gamma(\xi))
t}g(\lambda)d\lambda=\int_{0}^{+\infty}e^{-\lambda t}\tilde
g(\xi,\lambda)d\lambda,\ {\frak r}\ \textrm{a.e.}\ \xi,\label{koniec}
\end{align}
where
$$
\tilde g(\xi,\lambda)=\left\{ \begin{array}{cc}
\tilde\gamma(\xi)^{-1}g\left(\frac{\lambda-\gamma(\xi)}{\tilde\gamma(\xi)}\right),& \lambda>\gamma(\xi),\\
0,& \lambda\leq\gamma(\xi).
\end{array}\right.
$$
We can see that the last expression in \eqref{koniec} is of the form
\eqref{hfunkcja1} and the measure $\tilde\mu(\xi,d\lambda):=\tilde
g(\xi,\lambda)d\lambda$ fulfills \eqref{miaranos}.

\section{Proof of the Proposition \ref{dodatnioh}}\label{pomocniczedowody}
Denote by $A(\xi)$ some non-negative Hermitian, matrix valued
function such that
$$
\hat \bR(\xi)=A(\xi)\hat r(\xi),\quad \xi\in\mbR^d,
$$
where
\begin{equation}\label{slad}
{\rm tr}A(\xi)\equiv 1
\end{equation}
for ${\frak r}$ a.e. $\xi$. Choose a function $\hat\phi\in
\mathcal{S}(\mbR^d)$, where $\mathcal{S}(\mbR^d)$ is the class of
$\mathbb C^d$ valued Schwartz functions. Then
\begin{align*}
&0\le \bbE\left|\sum_{p=1}^N\al_p\int_{\mbR^d}\bbV(t_p,x)\cdot
\hat\phi(x) dx\right|^2=\frac{1}{(2\pi)^d}\int_{\mbR^d}
\left(\sum_{p,q=1}^N\al_p\bar\al_q h(t_p-t_q,\xi)\right)A(\xi)
\phi(\xi)\cdot\phi(\xi) \hat r(\xi)d\xi,
\end{align*}
where $\phi(\xi)$ is the inverse Fourier transform of a function
$\hat\phi(x)$. If we choose $\hat\phi_j(x)$, such that
$\phi_j(\xi)=e_j\psi(\xi)$, where $\psi(\xi)$ is a scalar function,
$e_j=(\underbrace{0,\dots,1,\dots,0}_{j\textrm{th position}})$ and
summing by $j=1,\dots,d$, we obtain
$$
0\leq\frac{1}{(2\pi)^d}\int_{\mbR^d}
\left(\sum_{p,q=1}^N\al_p\bar\al_q
h(t_p-t_q,\xi)\right)|\psi(\xi)|^2\underbrace{\sum_{j=1}^d A(\xi)
e_j\cdot e_j}_{=1} \hat r(\xi)d\xi.
$$
Considering \eqref{slad} and the fact that $\psi$ was arbitrary, we
obtain \eqref{pos-def}.\qed
\subsection*{Acknowledgement}
I would like to express my deep gratitude to Professor Tomasz
Komorowski whose help and comments have been invaluable.
 The work
has been supported by the NCN grant UMO 2016/23/B/ST1/00492.

\end{document}